\documentclass[11pt]{amsart}
\usepackage{amssymb}
\usepackage{microtype}
\usepackage[hidelinks]{hyperref}
\newtheorem{theorem}{Theorem}
\newtheorem{proposition}[theorem]{Proposition}
\newcommand{\F}{\mathbb F_3}

\begin{document}
\title[Bounds for extraspecial groups of exponent three]{Bounds for the small Davenport constant of extraspecial groups of exponent three}
\author{Andreas Volkmann}
\subjclass[2020]{11B30, 20D15}
\keywords{small Davenport constant, product-one subsequence, extraspecial group, subset sums}
\begin{abstract}
Let \(E_s\) be the extraspecial group of order \(3^{2s+1}\) and exponent three. We prove
\(4s+2r+2\le d(E_s\times C_3^r)\le 4s+2r+3\) for \(s\ge1\) and \(r\ge0\).
For \(s\ge2\), the degree in the polynomial proof is optimal among weights depending only on group-element multiplicities modulo three.
For \(s=2\), we reduce the question in all ranks to an abelian subset-sum problem with a fixed set of forbidden values in a \(77\)-dimensional space. We show that this set, together with zero, contains no \(76\)-dimensional subspace. The equality \(d(E_2\times C_3^r)=2r+10\) in all ranks is not established here.
\end{abstract}
\maketitle

\section{Definitions and results}
A sequence is an indexed list with repetitions allowed. A subset of its positions is product-one if its terms admit an ordering with product equal to the identity. We use \(d(G)\) for the largest length of a sequence over \(G\) with no nonempty product-one subset.

Write \(G_{s,r}=E_s\times C_3^r\) in coordinates \((v,h)\in V\times\F\), where \(V=\F^{2s+r}\) and \(h\) is the central coordinate corresponding to the commutator subgroup. The multiplication is
\[
 (v,h)(w,j)=(v+w,h+j+\tfrac12\omega(v,w)),
\]
where \(\omega\) is alternating of rank \(2s\) and has an \(r\)-dimensional radical. In this model \(\F\) is written additively. The case \(s=1\) has the stronger exact formula \(d(G_{1,r})=2r+6\) by \cite{volkmann-heis}. For \(s=2\), \cite{volkmann-small} proves the conjectured formula for \(0\le r\le3\), and \cite{volkmann-twelve} treats the layer with twelve noncentral terms at the critical length. The central compression in \cite{volkmann-central} reduces the rank by replacing a suitable central block. The construction in Section~\ref{sec:fixed} increases the rank and makes the entire noncentral sequence fixed.

\begin{theorem}\label{thm:bounds}
For \(s\ge1\) and \(r\ge0\),
\[
 4s+2r+2\le d(G_{s,r})\le 4s+2r+3.
\]
For each fixed \(s\), the function \(r\mapsto d(G_{s,r})-2r\) is nondecreasing. Hence the upper value, if attained at some rank, is attained at every larger rank.
\end{theorem}

\begin{proof}
Set \(m=2s+r\). Two copies of each of the \(m\) spatial basis vectors with height zero, together with two copies of \((0,1)\), give the lower bound. A product-one subset must omit all the spatial basis positions, and neither a single copy nor both copies of \((0,1)\) have product one.

Let \(S=((v_i,h_i))_{i=1}^n\) be product-one-free. On a spatially neutral subset \(J\), put an edge between \(i,j\in J\) when \(\omega(v_i,v_j)\ne0\). If there are two disjoint edges, arrange the two pairs adjacently. Switching their internal orders changes the height independently by two nonzero elements of \(\F\). Their four choices attain all three heights, which is impossible for a nonempty subset of \(S\). Each selected vertex satisfies \(\sum_{j\in J}\omega(v_i,v_j)=0\). Thus there is no vertex of degree one. A nonempty edge family in which all edges meet is contained in a star or is a triangle. The degree condition rules out the star and forces the three cyclic pairings of the triangle to be equal and nonzero. We call such a triangle balanced. Every remaining position is isolated. Its two possible product heights are the additive height \(\sum_{i\in J}h_i\) plus the two nonzero elements of \(\F\). Therefore a nonempty neutral subset is either commuting with nonzero additive height or a balanced triangle with isolates and additive height zero.

For Boolean position variables \(x_i\), write \(L_j(x)=\sum_i (v_i)_j x_i\), \(H(x)=\sum_i h_i x_i\), and let \(B(x)\) be the sum of \(x_i x_j x_k\) over all balanced triples. On \(\{0,1\}^n\),
\[
 \prod_{j=1}^m(1-L_j(x)^2)\bigl(1-H(x)^2-B(x)\bigr)
 =\prod_{i=1}^n(1-x_i).
\]
The first factor on the left detects spatial neutrality. On a nonempty neutral subset the preceding classification makes the second factor zero. Both sides are one at the empty subset. Multilinear reduction on the Boolean cube preserves the function and does not increase the degree. The right-hand side has unique multilinear degree \(n\), so \(n\le2m+3\).

To pass to \(r+1\), append two copies of a generator supported in a new direct central \(C_3\) factor. A product-one subset must omit both positions. This gives \(d(G_{s,r+1})\ge d(G_{s,r})+2\).
\end{proof}

Write \(\epsilon_{s,r}=d(G_{s,r})-(4s+2r+2)\). For \(r\ge1\), embed one radical direction of \(G_{s,r}\) as the first axis of an additional symplectic plane in \(G_{s+1,r-1}\). The form is unchanged on the embedded space. Adding two copies of the second axis preserves product-one-freeness, since their new coordinate forces both copies to be omitted in a product-one subset. Therefore \(\epsilon_{s,r}\le\epsilon_{s+1,r-1}\). In particular, equality at the lower bound for every group \(E_s\) would imply equality at the lower bound for every \(G_{s,r}\).

\section{An optimal degree for universal multiplicity weights}
For a sequence \(S\) let \(c(S)\in\F^{G_{s,r}}\) record the multiplicity of each group element modulo three. A universal weight is a polynomial \(P\) with \(P(0)=1\) such that \(P(c(S))=0\) whenever the entire nonempty sequence \(S\) admits no product-one ordering.

\begin{theorem}\label{thm:degree}
For \(s\ge2\), the least total degree of a universal weight on \(G_{s,r}\) is \(4s+2r+3\).
\end{theorem}

\begin{proof}
The left-hand side of the Boolean identity in the proof of Theorem~\ref{thm:bounds} gives the required weight in the multiplicity variables. The linear forms \(L_j\) and \(H\) depend only on multiplicities. A balanced triple uses distinct spatial values, so the number of positional triples of any such type is the product of the three type multiplicities. Its contribution to \(B\) is cubic. The classification in the preceding proof applies to any \(S\) whose entire nonempty sequence admits no product-one ordering. This proves the claimed upper bound on the degree.

For the lower bound retain only types \((v,0)\) for \(v\ne0\), with multiplicities \(x_v\), and the single type \((0,1)\), with multiplicity \(t\). Every ternary vector of these multiplicities has its representative sequence with multiplicities \(0,1,2\). Set \(\tau(x)=\sum_{v\ne0}x_vv\). A weight vanishes whenever \(\tau(x)\ne0\). Choose linear coordinates with the \(m\) entries of \(\tau(x)\) as the first \(m\) coordinates. In the unique reduced polynomial representation over \(\F\), the weight is the product of \(\prod_{j=1}^m(1-\tau_j(x)^2)\) and its restriction to \(\ker\tau\) with the coordinate \(t\). This follows from the description of the function on \(\tau(x)=0\) and its vanishing on the other fibers. Linear substitution and reduction using \(u^3=u\) do not increase degree. The two factors use disjoint coordinates, so their degrees add.

Suppose the degree were at most \(2m+2\). The restriction would then have the form
\(p_0(x)+t p_1(x)+a t^2\), with \(p_0\) quadratic and \(a\) constant. At \(x=0\) its values for \(t=0,1,2\) are \(1,0,0\), giving \(a=-1\). If \(x\) represents a commuting neutral sequence, both \(t=1,2\) give nonidentity total products, so \(p_0(x)=1\). If \(x\) represents a balanced triangle with isolates and spatial sum zero, then \(t=0\) gives no product-one ordering, so \(p_0(x)=0\). Consequently \(1-p_0\) would be a quadratic on \(\ker\tau\) taking value zero on all the indicated commuting profiles and one on all the indicated triangle profiles. Its value at zero is zero.

Restrict to types in any nondegenerate symplectic four-space of \(V\). The preceding value requirements persist. A quadratic on the spatial kernel of these \(80\) coordinates extends to the whole \(80\)-dimensional multiplicity space by composing it with a linear retraction onto that kernel. Its constant term remains zero. The certificate in Appendix~\ref{sec:certificate} annihilates every linear and quadratic monomial on these coordinates, while the weighted sum of the prescribed values is \(2\) in \(\F\). This is a contradiction.
\end{proof}

\section{A fixed noncentral sequence for \texorpdfstring{\(E_2\)}{E2}}\label{sec:fixed}
Let \(P=\F^4\) with its nondegenerate symplectic form and let \(A=\F^{P\setminus\{0\}}\). Write \(e_t\) for its standard basis. In \(G_{2,r+80}\) use coordinates \((t,\rho,a,h)\), where \(t\in P\), \(\rho\in\F^r\), \(a\in A\) and \(h\in\F\). Let \(T\) be the sequence consisting of two copies of \(p_t=(t,0,e_t,0)\) for every \(t\ne0\). The \(160\) terms of \(T\) are fixed.

\begin{proposition}\label{prop:normalize}
If \(S=((t_i,\rho_i,h_i))_{i=1}^n\) over \(G_{2,r}\) is product-one-free, then \(T\) together with the central sequence
\[
 Z=\bigl((0,\rho_i,-e_{t_i},h_i)\bigr)_{i=1}^n,\qquad e_0=0,
\]
is product-one-free over \(G_{2,r+80}\). The construction changes the length from \(n\) to \(n+160\).
\end{proposition}

\begin{proof}
Suppose an ordered nonempty product-one subset is chosen from \(T\) and \(Z\). Let \(J\) be the indices of its selected central terms, \(n_t\) the number of indices \(i\in J\) with \(t_i=t\), and \(s_t\in\{0,1,2\}\) the number of chosen \(p_t\) positions. The \(A\)-coordinates give \(s_t\equiv n_t\pmod3\) for each \(t\ne0\). In particular \(s_t\le n_t\), and \(J\) is nonempty because \(J=\varnothing\) forces all \(s_t=0\).

Replace the chosen \(p_t\) positions by \(s_t\) distinct original positions of type \(t\) from \(J\), retaining the order of the \(p_t\) positions. For each \(t\ne0\) group the remaining \(n_t-s_t\) original positions into triples of type \(t\). Each triple has zero symplectic sum and central product. Append these triples and the original positions of type zero from \(J\). The commutator correction in the resulting word is the correction of the chosen \(p_t\) positions. Its radical sum and additive height are exactly those of the selected central terms. Thus all of \(J\), with each position used once, has product one in \(G_{2,r}\).
\end{proof}

Put \(\tau(a)=\sum_{t\ne0}a_t t\) and let \(U=\ker(\tau)\oplus\F\), viewed in the center with \(\rho=0\). Then \(\dim U=77\) and the center modulo \(U\) has dimension \(k=r+4\). Let \(F\subset U\) be the set of negatives, in the additive center, of products of nonempty subsets of \(T\) whose products are central, with all orders allowed. This set is independent of \(r\). A subset of \(T\) is uniquely specified up to choices of copies by its multiplicity vector \(s\in\{0,1,2\}^{80}\). Its \(A\)-coordinate is \(s\) and its product is central exactly when \(\tau(s)=0\). Hence \(0\notin F\).

For a central sequence \(Z\), write \(\Sigma^*(Z)\) for its nonempty subset sums. Then \(T Z\) is product-one-free exactly when
\[
 \Sigma^*(Z)\cap(\{0\}\cup F)=\varnothing.
\]
Indeed a word using \(T\) must use a nonempty central subset to cancel its central product, since no nonempty subset of \(T\) has product one. A word using only \(Z\) has product one exactly at sum zero. It follows from Proposition~\ref{prop:normalize} that the all-rank conjecture for \(E_2\) is equivalent to the assertion that every central sequence of length \(2k+3\) meets \(\{0\}\cup F\) in its nonempty subset sums, for every \(k=r+4\ge4\). In the converse direction a central sequence avoiding \(\{0\}\cup F\), together with \(T\), is itself a counterexample in \(G_{2,r+80}\). This equivalence concerns existence of counterexamples across ranks. The proof of Proposition~\ref{prop:normalize} only requires the stated direction of transfer for an individual sequence.

Let \(a(k)\) be the maximum length of a sequence over \(U\oplus\F^k\) that avoids \(\{0\}\cup F\) in its nonempty subset sums, with the target set embedded in \(U\). This is the generalized Davenport constant of \cite{godinho} minus one. The two constructions, which preserve their stated lengths for any input sequence, give
\[
 d(G_{2,r})\le a(r+4)\le d(G_{2,r+80})-160 .
\]
In particular \(2k+2\le a(k)\le2k+3\) for \(k\ge4\). Appending a doubled basis vector in a new quotient direction gives \(a(k+1)\ge a(k)+2\). Both integer excesses are eventually constant, and the two bounds give
\[
 \lim_{r\to\infty}\bigl(d(G_{2,r})-2r\bigr)
 = 8+\lim_{k\to\infty}\bigl(a(k)-2k\bigr).
\]

\begin{proposition}\label{prop:subspace}
The set \(\{0\}\cup F\) contains no \(76\)-dimensional linear subspace of \(U\).
\end{proposition}

\begin{proof}
Over \(a=0\), the only point of \(\{0\}\cup F\) is \((0,0)\). Therefore any \(76\)-dimensional subspace contained in this set would be the graph of a linear form \(\ell\) on \(\ker\tau\). For a commuting multiplicity profile \(s\) with \(\tau(s)=0\), the only possible height over \(-s\) in \(\{0\}\cup F\) is zero. Thus \(\ell(-s)=0\) on all such profiles.

These profiles span \(\ker\tau\). To see this, consider an annihilator of their span, extend it to \(A\), and write \(f(t)\) for its value at \(e_t\), with \(f(0)=0\). The commuting profiles \(e_t+e_{-t}\) give \(f(-t)=-f(t)\). For \(\omega(t,u)=0\) and \(t+u\ne0\), the profile \(e_t+e_u-e_{t+u}\), interpreted with ternary multiplicities, gives \(f(t+u)=f(t)+f(u)\). The zero-sum and repeated-vector cases also follow from oddness. Such an orthogonally additive function is linear \cite[Lemma~1]{delfosse}. In dimension four one can see this directly by taking an orthogonal sum of two symplectic planes. Orthogonality makes \(f\) split as the sum of its restrictions. On either plane, the additivity defect for a nonorthogonal pair is determined by its pairing, by subtracting the defect of a pair with opposite pairing in the other plane. Simultaneously negating the pair preserves its pairing and negates the defect, so the defect is zero. Hence the annihilator of the commuting profiles consists exactly of pullbacks of linear forms on \(P\), proving the spanning assertion.

It follows that \(\ell=0\). Take a symplectic pair \(u,v\) with \(\omega(u,v)=1\) and \(s=e_u+e_v+e_{-u-v}\). It is a balanced triangle of zero spatial sum, with possible heights precisely the two nonzero elements of \(\F\). Thus \((-s,0)\notin\{0\}\cup F\), a contradiction.
\end{proof}

For comparison, a \(76\)-dimensional subspace contained in \(\{0\}\cup F\) would give the desired bound immediately by quotienting the center by it. The quotient would have rank \(k+1\), and \(d(C_3^{k+1})=2k+2\). This latter identity follows from the lower construction with doubled basis vectors and the augmentation ideal of \(\F[C_3^{k+1}]\), whose \((2k+3)\)-rd power vanishes. Proposition~\ref{prop:subspace} identifies a specific limitation of this quotient argument.

\section{Status of the all-rank problem}
Theorem~\ref{thm:bounds} leaves two consecutive values in each rank. For \(s=2\) the lower value is conjectured in every rank, and any occurrence of the upper value would persist in every larger rank. Proposition~\ref{prop:normalize} and the subset-sum formulation isolate the same unresolved choice as an abelian problem with fixed \(F\). The degree result concerns universal multiplicity weights. The subspace result concerns a particular quotient argument. Neither supplies a counterexample to the conjectured equality or a proof of it.

\appendix
\section{Certificate for the quadratic obstruction}\label{sec:certificate}
Use coordinates \((x_1,x_2,x_3,x_4)\) with the form \(x_1y_2-x_2y_1+x_3y_4-x_4y_3\). Each following row gives a weight, its prescribed value type (\(\mathrm B\) for one and \(\mathrm K\) for zero), and the multiset of spatial vectors. A superscript \(2\) indicates two copies of that vector. The coordinates of each vector are written consecutively.
\[
\begin{array}{rcl}
1&\mathrm B&0001,\ 0112,\ 0220\\
2&\mathrm B&0001,\ 0012,\ 0120,\ 0200\\
2&\mathrm B&0001,\ 0011,\ 0121,\ 0200\\
2&\mathrm B&0001,\ 0010,\ 0122,\ 0200\\
1&\mathrm B&0001,\ 0010,\ 0022\\
1&\mathrm B&0001,\ 0111,\ 0221\\
1&\mathrm B&0001,\ 0110,\ 0222\\
2&\mathrm B&0001,\ 0110,\ 1120,\ 2102\\
2&\mathrm B&0001,\ 0110,\ 1020,\ 2202\\
2&\mathrm B&0001,\ 0110,\ 1220,\ 2002\\
2&\mathrm K&2010,\ 2111,\ 2212\\
1&\mathrm K&2011,\ 2112,\ 2210\\
2&\mathrm K&1222,\ 2010,\ 2111^2,\ 2212\\
2&\mathrm K&1220,\ 2012,\ 2110^2,\ 2211\\
2&\mathrm K&1221,\ 2011,\ 2112^2,\ 2210
\end{array}
\]
For each row apply independently to the first and second pair of coordinates a matrix from the set consisting of
\[
 \begin{pmatrix}1&0\\0&1\end{pmatrix},\
 \begin{pmatrix}0&1\\2&0\end{pmatrix},\
 \begin{pmatrix}1&1\\1&2\end{pmatrix},\
 \begin{pmatrix}1&2\\2&2\end{pmatrix}
\]
and their negatives. Include both choices of whether to swap the resulting coordinate pairs. All \(128\) transformations preserve \(\omega\). If \(x\) is the multiplicity vector of a resulting row and \(w\) its inherited weight, the exact certificate identities in \(\F\) are
\[
 \sum w x_i=0,\qquad \sum w x_i x_j=0\quad(i\le j),\qquad
 \sum_{\mathrm B}w=2 .
\]
The supplied \texttt{degree\_certificate.py} checks the row types, all the monomial identities, and the last sum using integer arithmetic reduced modulo three.

\section*{Use of AI assistance}
The author used ChatGPT to develop and review arguments, assist with computational checks, and prepare the manuscript.

\end{document}